\documentclass[12pt]{article}

\usepackage{amsmath}
\usepackage{amssymb}
\usepackage[T1]{fontenc} 

\usepackage{colortbl}
\definecolor{text1}{cmyk}{1,.65,0,0} % blue text colour
\definecolor{text2}{rgb}{1,0,0} % red text colour
\definecolor{text3}{cmyk}{0,0,0,1} % black text colour
\definecolor{text4}{cmyk}{0,0,0,0.5} % grey text colour

\usepackage{amsfonts}
\usepackage{amsthm}
\usepackage{graphicx}
\usepackage{latexsym}
\makeatletter
\renewcommand{\@seccntformat}[1]
{\csname the#1\endcsname.\enspace}
\makeatother
\newtheorem{theorem}{Theorem}
\newtheorem{lemma}{Lemma}
\newtheorem{remark}{Remark}
\newtheorem{corollary}{Corollary}
\newtheorem{example}{Example}

\begin{document}

\begin{center}
   {\bf Inference for a constrained parameter in presence of an uncertain constraint \footnote{ \today}}
\end{center}

\begin{center}
{\sc  \'Eric Marchand$^{a}$, Theodoros Nicoleris$^{b}$}\\

{\it a  Universit\'e de
    Sherbrooke, D\'epartement de math\'ematiques, Sherbrooke Qc,
    CANADA, J1K 2R1} \quad 
    
{\it b  Department of Economics, National and  Kapodistrian  University of Athens,
GR-105 59     Athens, GREECE \quad     
    (e-mails: eric.marchand@usherbrooke.ca; tnicoleris@econ.uoa.gr) } \\

\end{center}
\vspace*{0.2cm}
\begin{center}
{\sc Summary} \\
\end{center}
\vspace*{0.1cm}
\small
We describe a hierarchical Bayesian approach for inference about a parameter $\theta$ lower-bounded by $\alpha$ with uncertain $\alpha$, derive some basic identities for posterior analysis about $(\theta,\alpha)$, and provide illustrations for normal and Poisson models.  For the normal case with unknown mean $\theta$ and known variance $\sigma^2$, we obtain Bayes estimators of $\theta$   
that take values on $\mathbb{R}$, but that are equally adapted to a lower-bound constraint in being minimax under squared error loss for the constrained problem.

\vspace*{0.5cm}


\section{Introduction}

Consider a statistical model $X \sim  f_{\theta}$ and the problem of estimating $\theta$ under the parametric constraint $\theta \in C \subset \mathbb{R}^p$.   There are many challenging aspects to making inferences in such restricted parameter space settings and an accompanying rather large literature (e.g., Mandelkern, 2002; Marchand and Strawderman, 2004; van Eeden, 2006) on the frequentist performance of various estimators.  But they deal mostly with settings for fixed $C$; where there is no uncertainty in the parametric information and the objective is how to capitalize on such information for inferential purposes.      
This notes deals with situations where uncertainty resides in the parametric constraint and where we take a simple hierarchical Bayes approach to describe the uncertainty relative to $C$ and to $\theta$.  

\noindent More specifically, for estimating the mean $\theta$ of a normal distribution $N(\theta, \sigma^2)$ with known $\sigma^2$, under squared-error loss,  we obtain a class of hierarchical Bayes estimators that dominate $X$ on the restricted parameter space $[0,\infty)$.  The concerned priors are of the form:
\begin{equation}
\label{priorsthatwork}
\pi(\theta|\alpha) \, = \, \mathbb{I}_{[\alpha, \infty)}(\theta)\,,\, \alpha \sim N(0, \sigma^2_{\alpha})\,, \hbox{ with } \sigma^2_{\alpha} >0\,.
\end{equation}
The resulting Bayes point estimators of $\theta$ take values on $\mathbb{R}$, but are also adapted to performing better than the unbiased estimator $X$ under a lower bound constraint $\theta \geq 0$, achieving strict dominance and thus being minimax given that $X$ is minimax.  The limiting case $\sigma^2_{\alpha} \to 0$ yields Katz (1961) estimator shown by Katz to be minimax on $[0,\infty)$, while the limiting case $\sigma^2_{\alpha} \to \infty$ leads to $X$.  The finding is obtained via Stein's identity, sign change arguments and properties of the inverse Mill's ratio which intervenes in the functional form of the Bayes estimators.  We believe the finding is significant as such an interaction between the Bayes character and the minimax property in presence of a lower bound is particularly interesting and seems to have been undiscovered up to now.  We do point out that it has long been recognized (e.g., O'Hagan and Leonard, 1975) that placing a hierarchical prior on the bounds that may apply to an unknown parameter is an attractive choice for proceeding with Bayes inference for such problems with uncertain constraints. 

The paper is organised as follows.  In Section 2, we present a simple hierarchical framework to reflect uncertainty in a pre-specified lower bound constraint, extract some basic properties and present some examples which are of independent interest, including the Poisson and normal cases with various skewed-type distributions arising. 
  In Section 3, we establish and comment on the main minimaxity result for a class of hierarchical Bayes estimators of a lower-bounded normal mean.  Finally, Section 4 contains some final observations.

\section{Priors, posterior analysis and examples}
 
\noindent We focus on the case of a lower bound constraint with model and prior represented as:
%\begin{equation}
%\label{model0}
%X|\theta, \alpha \sim f_{\theta}(x)\,; \;
%(\theta,\alpha) \sim g(\theta,\alpha) \,\mathbb{I}_{[0,\infty)}(\theta-\alpha)\,.
%\end{equation}
\begin{equation}
\label{model}
X|\theta, \alpha \sim f_{\theta}(x)\,; \theta|\alpha \sim g_1(\theta) \,  \mathbb{I}_{[0,\infty)}(\theta-\alpha)\, , \, \alpha \sim g_2(\alpha)\,,
\end{equation}

The above describes indeed a situation where $\theta$ is lower bounded by  
 $\alpha$, but $\alpha$ is viewed as random or uncertain.  It also describes the conditional distribution of $X|\theta, \alpha $ as independent of
 $\alpha$.   Cases where the marginal prior distribution of $\alpha$ is degenerate reduce to cases where the lower bound constraint is deterministic.   
%  We will also simplify the prior in (\ref{model0}) to the multiplicative case :
We further assume that $g_1$ is absolutely continuous with respect to Lebesgue density and that $g_2$ is absolutely continuous with respect to a $\sigma-$finite measure $\mu$ with finite $\int_{-\infty}^{t} \, g_2(\alpha) \, d\mu(\alpha) < \infty$, for all $t \in \mathbb{R}$.  

\begin{remark}
\label{priordensityequivalent}
If the focus is inference on $\theta$, which is the case in Section 3, one can take equivalently the prior 
\begin{equation}
\label{priorequivalent} 
 \pi(\theta) \, = \, g_1(\theta) \, \int_{-\infty}^{\theta} \, g_2(\alpha) \, d\mu(\alpha)\,,
\end{equation} 
reducing to the skewed version $g_1(\theta)\, G_2(\theta)$ of $g_1$ in cases where $g_2$ is a density with c.d.f. $G_2$.  Alternatively, inferential interest may lie with $\alpha$, but this will at least require the finiteness of $\int_{s}^{\infty} \, g_1(\theta) \, d\theta $ for all $s \in \mathbb{R}$.  We do not necessarily assume this, allowing for instance the choice $g_1(\theta)=1$, but we will implicitly assume this finiteness when describing the posterior distribution of $\alpha$ (such as in Theorem \ref{posterior}).  \end{remark}

The above model leads to the following.
\begin{theorem}  
\label{mainposterior}
Under a model and prior as in (\ref{model}), the marginal posterior densities of $\theta$ and $\alpha$ are (whenever they exist) given respectively as 
\begin{eqnarray}
\label{posttheta}
\pi_1(\theta|x) &\propto&   f_{\theta}(x) \, g_1(\theta) \, 
\int_{-\infty}^{\theta} \, g_2(\alpha) \, d\mu(\alpha)\,,\\
\label{postalpha}
\pi_2(\alpha|x) &\propto&  g_2(\alpha)  \int_{\alpha}^{\infty} \, 
f_{\theta}(x) \, g_1(\theta)\,d\theta \,.
\end{eqnarray}
\end{theorem}
\noindent {\bf Proof.}   These expressions follow immediately by expressing the joint posterior as $\pi(\theta,\alpha) \propto f_{\theta}(x) \, g_1(\theta) \, g_2(\alpha) \, \mathbb{I}_{[0,\infty)}(\theta-\alpha)$.  \qed

\begin{remark}
%\begin{enumerate}
%\item[ \bf{(a)}]
Let $\pi_0(\theta|x)$ be the posterior density of $\theta$ in the absence of a lower bound constraint for the prior (i.e., $\alpha=-\infty$).  The posterior density in (\ref{posttheta}) is a weighed version of $\pi_0$ expressible as
\begin{equation*}
\label{priortheta}
\pi_1(\theta|x)\,=\, \pi_0(\theta|x) \, w(\theta),
\end{equation*}
with weight $w(\theta) = \int_{-\infty}^{\theta} \, g_2(\alpha) \, d\mu(\alpha)\,.$    Observe that $w(\cdot)$ is non-decreasing so that the posterior distributions $\pi_1(\theta|x)$ and $\pi_0(\theta|x)$ are stochastically ordered in terms of a m.l.r. with $\pi_1(\cdot|x)$ stochastically larger.  

\end{remark}

\subsection{Example (normal model)}

\subsubsection{Estimation of $\theta$}

Consider the normal case in (\ref{model}) with 
$X|\theta,\alpha \sim N(\theta,\sigma^2)$, known $\sigma^2$, denoting  $\phi$ and $\Phi$ as the probability density and cumulative distribution functions of a $N(0,1)$ random variable.
Furthermore, consider a normally distributed $g_1 \sim N(\mu,\tau^2)$, so that $\theta|\alpha$ is distributed as truncated normal on $[\alpha,\infty)$, and $\alpha \sim N(0,\sigma_{\alpha}^2)$ with 
$\sigma_{\alpha}^2 \geq 0$; the degenerate case $\sigma_{\alpha}^2=0$ covering the deterministic non-negativity constraint $\theta \geq 0$.  Without the constraint, the posterior density is equal to 
\begin{equation}
\label{pi0}
\pi_0(\theta|x) \sim N(\hat{\mu}(x), \tau'^2)\,,
\end{equation}
with 
\begin{equation}
\label{parameters}
\hat{\mu}(x) = \frac{\tau^2}{\tau^2 + \sigma^2} \, x + \frac{\sigma^2}{\tau^2 + \sigma^2}\, \mu \hbox{ and }  \tau'^2=\frac{\sigma^2\,\tau^2}{\tau^2 + \sigma^2}\,.  
\end{equation}
With the above choice of prior on $(\theta,\alpha)$, we obtain from (\ref{posttheta})
\begin{eqnarray}
\label{post7} \pi_1(\theta|x) \propto \phi(\frac{\theta - \hat{\mu}(x)}{\tau'}) \; \Phi(\frac{\theta}{\sigma_{\alpha}})\,\,\, \hbox{for } \sigma^2_{\alpha} >0\,, \\
\nonumber
\hbox{and } \pi_1(\theta|x) \propto \phi(\frac{\theta - \hat{\mu}(x)}{\tau'}) \,\mathbb{I}_{[0,\infty)}(\theta) \;\,\,\, \hbox{for } \sigma^2_{\alpha} =0\,.
\end{eqnarray} 
The following describes further the densities in (\ref{post7}).  These densities and given properties given are familiar (see Remark \ref{normalization}), but collected here for completeness. 

\begin{theorem}
\label{posterior}
The posterior density of $Z=\frac{\theta-\hat{\mu}(x)}{\tau'}$ is equal to 
\begin{equation}
\label{posteriordensityZ}
 f_{\psi_1,\psi_2}(z)\,=\,\frac{\phi(z) \, \Phi(\psi_1+\psi_2 z)}{\Phi(\frac{\psi_1}{\sqrt{1+ \psi_2^2}})}\,,
\end{equation}
with $\psi_1=\frac{\hat{\mu}(x)}{\sigma_{\alpha}}$ and $\psi_2=\frac{\tau'}{\sigma_{\alpha}}$.  Furthermore, we have, denoting $R(t)=\frac{\phi(t)}{\Phi(t)}, t \in \mathbb{R},$ the inverse Mill's ratio,
\begin{equation}
E[e^{tZ}]\,=\, e^{\frac{t^2}{2}} \frac{\Phi(\gamma_1 t + \gamma_0)}{\Phi(\gamma_0)} \,,\; E(Z) = \gamma_1 \, R(\gamma_0)\,, \; 
\hbox{and  Var}(Z) \, = \, 1 -  {\gamma_1}^2 \, R(\gamma_0) \, (\gamma_0 + R(\gamma_0))\,,
\end{equation}
$t \in \mathbb{R}$, with $\gamma_0= \frac{\psi_1}{\sqrt{1+\psi_2^2}}$ and $\gamma_1 = \frac{\psi_2}{\sqrt{1+\psi_2^2}}$. 
\end{theorem}
{\bf Proof.}
The given expressions for
$E(Z)$ and $\hbox{Var}(Z)$ follow readily by taking two derivatives of the moment generating function, while the normalization constant in (\ref{posteriordensityZ}) follows by taking $t=0$ in the development below.  For the moment generating function, we have 
\begin{eqnarray*}
E[e^{tZ}]\,&\propto\,&  \, \int_{\mathbb{R}} e^{tz} \phi(z) \int_{-\infty}^{\psi_1+\psi_2 z} \phi(w) \,dw \,dz  \\
\,& \propto \,& e^{\frac{t^2}{2}} \, \int_{\mathbb{R}} \phi(z-t) \, \int_{-\infty}^{\psi_1+\psi_2z} \phi(w)  \,dw \,dz\\
%\,&=\,& \frac{e^{\frac{t^2}{2}}}{\Phi(\gamma_0)}  \int_{\mathbb{R}}  \phi(w) \, \Phi(t+\frac{a-w}{b})  \, dw \,,
\,& \propto \,& e^{\frac{t^2}{2}} \;  \mathbb{P}(U_2 - \psi_2 U_1 \leq \psi_1)   \,,
\end{eqnarray*}
with $U_1 -t, U_2$ independent $N(0,1)$; and the result follows.  \qed
\begin{remark}
The developments above are also applicable to the case where the (improper) density of $\theta|\alpha$ is constant on $[\alpha,\infty)$ by taking $\hat{\mu}(x)=x$ and $\tau'^2=1$.  This may be viewed by taking $\tau^2 \to \infty$ for the density $g_1$.
\end{remark}

\begin{remark}
\label{normalization}
The densities in (\ref{posteriordensityZ}) coincide with a class of densities introduced by Azzalini (1985) and further analyzed by
Arnold et al. (1993), with the particular case $\psi_1=0$ reducing to the original skew normal density 
$2\phi(z) \, \Phi(\psi_2 z)$ introduced in Azzalini's seminal 1985 paper.  It is particularly interesting that such skewed-normal distributions arise in our setting.  Actually, the results of this section can be summarized as follows.  For $\theta|\alpha \sim N (\mu, \tau^2)$ truncated to $[\alpha,\infty)$, $\alpha \sim N(0, \sigma^2_{\alpha})$,  we have from (\ref{priorequivalent}) as a prior $$Z'= \frac{\theta-\mu}{\tau} \sim f_{\frac{\mu}{\sigma_{\alpha}}, \frac{\tau}{\sigma_{\alpha}}}\,,$$ 
and, from Theorem \ref{posterior}, as a posterior 
$$Z = \frac{\theta - \hat{\mu}(x)}{\tau'} \, \sim f_{_{\frac{\hat{\mu}(x)}{\sigma_{\alpha}}, \frac{\tau'}{\sigma_{\alpha}}}}\,.$$

%The normalization constant may be derived by the familiar expansion $\int_{\mathbb{R}} \phi(z) \Phi(\psi_1 + \psi_2 z) \, dz = \int_{\mathbb{R}} \phi(z) \int_{-\infty}^{\psi_1 + \psi_2 z} \phi(y) dy \, dz = P(Y_1-\psi_2 Y_2 \leq \psi_1 )=
%\Phi(\frac{a}{\sqrt{1+ b^2}})$, with $Y_1, Y_2$ independent N$(0,1)$.
\end{remark}

\noindent  From Theorem \ref{posterior}, we obtain the following Bayes estimators.

\begin{corollary}
\label{bayesestimator}
For $X|\theta \sim N(\theta, \sigma^2)$, $g_1 \sim N(\mu, \tau^2)$ and 
$\alpha \sim N(0,\sigma_{\alpha}^2)$, the Bayes point estimator of $\theta$ under loss $(d-\theta)^2$ is given by
\begin{equation}
E(\theta|x) \,=\, \hat{\mu}(x) \, + \, \frac{\tau'^2}{\sqrt{\tau'^2 + \sigma_{\alpha}^2}} \, R(\frac{\hat{\mu}(x)}{\sqrt{\tau'^2 + \sigma_{\alpha}^2}})\,
\end{equation}
with $\hat{\mu}(x)$ and $\tau'^2$ given in (\ref{parameters}).
%\footnote{For the case $\alpha \sim N(\mu_{\alpha},\sigma_{\alpha}^2)$, we obtain
%\begin{equation}E(\theta|x) \,=\, \hat{\mu}(x) \, + \, \frac{\tau'^2}{\sqrt{\tau'^2 + \sigma_{\alpha}^2}} \, R(\frac{\hat{\mu}(x)- \mu_{\alpha}}{\sqrt{\tau'^2 + \sigma_{\alpha}^2}})\,.
%\end{equation}} 
For the case of the uniform prior $\pi(\theta|\alpha)\,=\,\mathbb{I}_{[\alpha,\infty)}(\theta)$, the Bayes point estimator is as above with $\hat{\mu}(x)=x$ and $\tau'^2=\sigma^2$, that is
\begin{equation}
\label{bayesflat}
E(\theta|x) \,=\, x \,+\, \frac{\sigma^2}{\sqrt{\sigma^2 + \sigma_{\alpha}^2}} \, R(\frac{x}{\sqrt{\sigma^2 + \sigma_{\alpha}^2}})\,.
\end{equation}
\end{corollary}
 
\subsubsection{Estimation of the lower bound $\alpha$}

Interest may reside in estimating the lower bound as well.  As an illustration, consider again (\ref{model}) with 

\begin{equation}
\label{set-up}
X|\theta,\alpha \sim N(\theta, \sigma^2)\,,\, g_1 \sim N(\mu, \tau^2), \,\,
g_2 \sim N(\mu_{\alpha}, \sigma_{\alpha}^2)\,. 
\end{equation}

\begin{corollary}
For the normal model (\ref{set-up}-\ref{model}),
the posterior density of $W=\frac{\alpha - \mu_{\alpha} }{\sigma_{\alpha}}$
is given by $f_{\psi_1,\psi_2}$ in (\ref{posteriordensityZ}) with $\psi_1=\frac{\mu_{\alpha}-\hat{\mu}(x)}{\tau'}$ and $\psi_2=\frac{\sigma_{\alpha}}{\tau'}$, and $\hat{\mu}(x)\,,\, \tau'$ given in (\ref{parameters}).   The posterior expectation is equal to 
$$ E(\alpha|x) \,=\,  \mu_{\alpha} - \frac{\sigma_{\alpha}^2}{\sqrt{\tau'^2 + \sigma_{\alpha}^2}} \, R(\frac{\mu_{\alpha}-\hat{\mu}(x)}{\sqrt{\tau'^2 + \sigma_{\alpha}^2}}) \,.$$  
\end{corollary}
\noindent {\bf Proof.}  The posterior expectation follows from Theorem \ref{posterior}.  The posterior density of $\alpha$ may be derived from (\ref{postalpha}) with
$$\pi_2(\alpha|x) \propto \phi(\frac{\alpha - \mu_{\alpha}}{\sigma_{\alpha}}) \, \int_{\alpha}^{\infty} \, \pi_0(\theta|x) \,d\theta\,,$$
and $\pi_0$ given in (\ref{pi0}).  The posterior distributions for 
both $\alpha$ and $W$ follow.  \qed

\subsection{Example (Poisson model)}

In (\ref{model}), we consider the Poisson case with $X|\theta,\alpha \sim \hbox{Poisson}(\theta)$ and present analysis and inference with the choices
(which include Gamma densities)
\begin{equation}
\label{priorpoisson}
  g_1(\theta) \propto \theta^{a-1} \, e^{-b\theta} \, \mathbb{I}_{[0,\infty)}(\theta)\;,\; g_2(\alpha) \propto \alpha^{c-1} \, e^{-d\alpha} \, \mathbb{I}_{[0,\infty)}(\alpha)\,,
\end{equation}
with $a,c>0$, $b > -1$.  Several interesting cases arise.  Heuristically, a sequence of choices for $g_2$ such that $d \to \infty$ will lead to the unrestricted case with $\alpha=0$ with probability one.  
Denote  $f_{\gamma_1, \gamma_2}, F_{\gamma_1, \gamma_2},$ and 
$\bar{F}_{\gamma_1, \gamma_2}$ as the probability density, cumulative distribution and survivor functions (respectively) of a Gamma$(\gamma_1, \gamma_2)$ distribution.  The following corollaries are consequences of Theorem \ref{mainposterior}.
\begin{corollary}
\label{poissonposteriors}
For model (\ref{model}) with $X|\theta,\alpha \sim \hbox{ Poisson}(\theta)$,
$g_1$ and $g_2$ as in (\ref{priorpoisson}), we have
\begin{eqnarray*}
\pi_1(\theta|x) \propto \theta^{a+x-1} \, e^{-\theta(1+b)} \, \int_0^{\theta} \, \alpha^{c-1} \, e^{-d\alpha} \, d\alpha\,, \\
\hbox{or } \pi_1(\theta|x) \propto  \theta^{a+x-1} \, e^{-\theta(1+b)} \, F_{c,d}(\theta)\,,
\end{eqnarray*}
whenever $d>0$.  
\end{corollary}
We point out that the above density a weighted version of a Gamma$(a+x,1+b)$ density, which is itself  recovered by taking $d \to \infty$.  Some cases lead to closed forms, namely: (i) $d=0$, and (ii) $d>0$ with integer $c$.  For instance, the case $c=1$ yields $\pi_1(\theta|x) = k \theta^{a+x-1} \, e^{-\theta(1+b)} \, (1-e^{-d\theta})$ with $k^{-1}= \Gamma(a+x) \, 
\left(\frac{1}{(1+b)^{a+x}} - \frac{1}{(1+b+d)^{a+x}}  \right)$.        

\begin{corollary}   
For model (\ref{model}) with $X|\theta,\alpha \sim \hbox{ Poisson}(\theta)$,
$g_1$ and $g_2$ as in (\ref{priorpoisson}), we have
\begin{equation}
\pi_2(\alpha|x)  \propto \alpha^{c-1} \, e^{-d\alpha} \, \bar{F}_{x+a,1+b}(\alpha)\,.
\end{equation}
Furthermore, when $a$ is a positive integer, we have the finite Gamma mixture representation:
\begin{equation}
\pi_2(\alpha|x) \,=\, \sum_{y=0}^{x+a-1} p_y \, f_{c+y,1+b+d}(\alpha)\,,
\end{equation}
with $p_y \propto \, \frac{\rho^y\, \Gamma(c+y)}{y!} \, \mathbb{I}_{\{0, \ldots, x+a-1\}}(y)\,$, $\rho= \frac{1+b}{1+b+d}\,$.
\end{corollary}
\noindent {\bf Proof.}   The results follow easily by applying Theorem \ref{mainposterior} and by making use of the closed form representation 
$\bar{F}_{x+a,1+b}(\alpha) \, = \, \sum_{y=0}^{x+a-1} \, 
\frac{e^{-\alpha(1+b)} \, (\alpha(1+b))^y }{y!}$.  \qed

\begin{remark}
\label{negbinomial}

\begin{enumerate}
\item[ (a)]
For $d>0$ and positive integer $a$, the mixing proportions $p_y$ are those of a truncated (to $\{0, \ldots, x+a-1\}$) Negative Binomial distribution and expressible as $p_y = P(Y=y|Y \leq x+a-1)$ with $P(Y=y)\,=\, \frac{\Gamma(c+y)}{y! \, \Gamma(c)} \, \rho^y \, 
(1-\rho)^c \, \mathbb{I}_{\mathbb{N}}(y)\,$.  From the representation, the posterior expectation of $\alpha$ is equal to
$$  E(\alpha|x) \,=\, \sum_y \, p_y \frac{c+y}{1+b+d}  = 
\frac{c}{1+b+d} \, + \, \frac{1}{1+b+d} \, E(Y|Y \leq a+x-1\,)\,,$$
bringing into play the moments of a truncated Negative Binomial distribution.
In the very specific case where we observe $X=0$ and where $a=1$, the above finite mixture is degenerate and the posterior distribution for $\alpha$ reduces to a single Gamma$(c,1+b+d)$ distribution.

\item[ (b)]  For the particular case of an entirely flat prior for $\alpha$ (i.e., $c=1$, $d=0$), the mixture proportions $p_y$ are simply those of a uniform distribution on $\{0, \ldots, x+a-1\}$ and we obtain easily, for instance, $E(\alpha|x) \,=\,\frac{1}{2(1+b+d)} \, (2c+x+a-1)$.  
\end{enumerate}
\end{remark}

We conclude this section by describing a potential application.
\begin{remark}
Setting $g_2$ as a mixture of a mass at $0$ with an absolutely continuous part on $(0,\infty)$ will yield a posterior which is also such a mixture, with the posterior probability $P(\alpha=0|x)$  helping to gauge the probability that there exists a constraint.   With respect to estimating $\theta$, such a prior does not place mass zero on any interval subset of the $\mathbb{R}_+$ which is potentially appealing and in contrast to a deterministic lower-bound constraint.
\end{remark}

\section{Minimaxity of a class of hierarchical Bayes estimators of 
$\theta$ under the restriction $\theta \geq 0$}

\noindent  Notice that the Bayes estimator with respect to the flat prior 
on $[0,\infty)$ for $\theta$ ($\alpha$ degenerate at $0$) is recovered as
$\delta_U(x)=x+ \sigma R(\frac{x}{\sigma})$ by taking $\sigma_{\alpha}^2=0$ above.  This estimator was considered by Katz (1961) who showed that it is admissible, as well as minimax, with minimax risk given by $\sigma^2$.  Moreover, both the maximum likelihood estimator $\max(0,X)$ and the unbiased estimator $\delta_0(X)=X$ are  minimax, although the former dominates the latter, (see Marchand and Strawderman, 2012, for instance for further details on such a phenomenon with varying models and losses).  

Now, consider of estimators of the form $\delta_c(x)= x+ c\sigma R(c x/\sigma)$ with $c \in [0,1]$, which include the unbiased estimator ($c=0$), the estimator 
$\delta_U$ ($c=1$), as well as all the Bayes estimators in (\ref{bayesflat}) for $c=\frac{\sigma}{\sqrt{\sigma^2+\sigma_{\alpha}^2}}$).  With $\delta_c$ minimax for $c=0$ and $c=1$ for the restriction $\theta \geq 0$, we show below that,  for all $c \in (0,1)$, $\delta_c$ dominates $\delta_0$ and is thus minimax  for the restricted parameter space $\theta \geq 0$.  This is particularly interesting result since $\delta_c(X)$: {\bf (I)} is a Bayes estimator taking values everywhere on $\mathbb{R}$, {\bf (II)} yet is adapted to the existence of an uncertain lower bound constraint on $\theta$, {\bf (III)} and performs as a minimax estimator dominating $\delta_0$ for $\theta \geq 0$.  In fact, it also dominates $\delta_0$ for $\theta \geq \theta_0(c)$ for some $\theta_0(c) <0$ and performs better than $\delta_U$ on and near the boundary of the parameter space.  This is true for all $c \in (0,1)$.

The following theorem is the main result of this section.  We make use of the following well-known properties of the inverse Mill's ratio.

\begin{lemma}
\label{inversemill}
The inverse Mill's ratio $R \equiv \frac{\phi}{\Phi}$ is a nonincreasing and convex function on $\mathbb{R}$.  Furthermore, $\lim_{t \to \infty} R(t)=0$, $R(t) \geq -t$ for
all $t \in \mathbb{R}$, $\lim_{t \to -\infty} \frac{R(t)}{t}=-1$,  $R'(t) = -R(t)(t+R(t))$, $\lim_{t \to -\infty} -R(t)\,(t+R(t))\,=\,-1$, and
$\lim_{t \to +\infty} R(t)\,(t+R(t))\,=\,0$.
\end{lemma}

\begin{theorem}
\label{minimaxdeltac}
For $X \sim N(\theta,\sigma^2)$, loss $(d-\theta)^2$, parameter space $\theta \geq 0$,
the estimators $\delta_c(X)=X+c\sigma R(cX/\sigma)$, $c \in (0,1]$, dominate $\delta_0(X)=X$ and are thus minimax.
\end{theorem}
\noindent {\bf Proof.}  Since $\delta_0(X)=X$ is minimax (e.g., Katz, 1961), it suffices to show that $\delta_c(X)$ dominates $X$ for $c \in (0,1]$.   Since 
$$ \mathbb{E} \left((X+ c\sigma R(\frac{cX}{\sigma}) - \theta)^2  \right) \, = \,  
\sigma^2 \, \mathbb{E} \left((Z+ c R(cZ) - \frac{\theta}{\sigma})^2  \right)\,,  $$
with  $Z \sim N(\theta_Z=\frac{\theta}{\sigma},1)$ and $\theta_Z \geq 0$, we can take $\sigma^2=1$ without loss of generality.   We proceed below to show that
\begin{enumerate}
\item[ {\bf (i)}]  $\Delta_c(\theta)=R(\theta,\delta_c) - R(\theta,\delta_0)$ changes signs at most once from $+$ to $-$ as $\theta$ varies on 
$\mathbb{R}$;
\item[ {\bf (ii)}]  $\Delta_c(0) \leq 0$ for all $c \in (0,1]$;
\end{enumerate}
which taken jointly will imply the result.  For {\bf (i)}, we first apply Stein's integration by parts identity (i.e., $E_{\theta}[h(X) (X-\theta)]=E_{\theta}[h'(X)]$ for (weakly) differentiable $h$ and subject to the existence of both expectations) to obtain
\begin{eqnarray}
\nonumber  \Delta_c(\theta) &=& E_{\theta}[(X+cR(cX)-\theta)^2] \, - \,
E_{\theta}[(X-\theta)^2] \\
\nonumber \, &=&  E_{\theta}[ c^2 R^2(cX) + 2 \frac{\partial }{\partial X}
(c R(cX)] \\
%\nonumber &=& E_{\theta}[ -c^2 R(cX) \left( R(cX) + 2cX 
%\right)] \\
\label{Deltac} &=& - c^2 \, E_{\theta}[T(cX)]\,, \, \hbox{ with } T(s)=R(s) \,(R(s)+2s).  
\end{eqnarray}
Now observe that $R(s) +2s$ increases in $s \in \mathbb{R}$, with limits of
$\pm \infty$ when $s \to \pm \infty$ respectively (Lemma \ref{inversemill}), which implies that 
$-c^2 T(cs)$ changes signs once from $+$ to $-$ as $s$ varies from $-\infty$ to $+\infty$.   
With a normal model for $X$ and with the possible changes of $\Delta_c(\theta)$ as a function of $\theta \in \mathbb{R}$ governed by the sign changes of $-c^2T(cs)$ (i.e., Karlin, 1957; Brown, Johnstone and MacGibbon, 1981), we infer that $\Delta_c(\theta)$ changes signs at most once, from 
$+$ to $-$, as a function of $\theta \in \mathbb{R}$ (and also for $\theta \in \mathbb{R}_+$), establishing {\bf (i)}.

For part {\bf (ii)}, we have from (\ref{Deltac})
\begin{eqnarray*}
-\frac{1}{c^2} \, \frac{\partial}{\partial c} \, \Delta_c(0) \,&=&\, c E_0[X T'(cX)] \\
\, &=& c \, \int_{\mathbb{R}_+} x \phi(x) \, \left(T'(cx)-T'(-cx) \right) \, dx \,.
\end{eqnarray*}
Since $\Delta_1(0)=0$ (e.g., Marchand and Strawderman, 2005),  {\bf (ii)} will follow if we can show that
\begin{equation}
\label{T'}  T'(s) \leq T'(-s)  \; \hbox{ for all } s >0\,.
\end{equation}
With $T'(s) = 2 R(s) \, \{1- (s+R(s))^2\}$, (\ref{T'}) is equivalent to
\begin{equation}
\label{T'2}
R(s) \{ (s+R(s))^2 \,-\,1\} \geq R(-s) \{ (-s+R(-s))^2 \,-\,1\} \; \hbox{ for all } s >0\,.
\end{equation}
Notice that the right-hand side of (\ref{T'2}) is negative since $(x+R(x))^2$ increases in $x$ and consequently $(-s + R(-s))^2 \leq (0 + R(0))^2 = \frac{2}{\pi} < 1$ for all $s>0$.   This, by virtue of the monotone increasing property of 
$(x+R(x))^2$, implies inequality (\ref{T'2}) for all $s \geq s_0$ where $s_0+R(s_0)=1$.  Finally, (\ref{T'2}) holds for all $s$, since $0 < R(s) \leq R(-s)$ and 
$0< 1- (s+R(s))^2 \leq 1- (-s+R(-s))^2$   for all $s \in (0,s_0)$.  \qed \\

\begin{example}
As an illustration, Figure 1 represents the risks of $\delta_c$ for $c=1/2,3/4,1$ as functions of 
$\theta$ for $X \sim N(\theta,1)$.  The choice $\delta_1$ (yellow) is Katz's minimax estimator for the restriction $\theta \geq 0$ with no uncertainty on the lower bound $\alpha$.  The minimax risk is $1$ and the choices $\delta_{1/2}$ (green) and $\delta_{3/4}$ (red) are also minimax for $\theta \geq 0,$ by virtue of Theorem \ref{minimaxdeltac}, with strict dominance as well.   These estimators are more robust when it turns out that $\theta <0$ than $\delta_1$, with $\delta_{1/2}$ the more robust of the two.   The latter still improves on $X$ for $\theta \in (\theta_0,0)$ with $\theta_0 \approx -0.939$ (a little less than one standard deviation away from $0$).   On the other hand, the gains offered by $\delta_{1/2}$ are less pronounced for $\theta \geq 0$.  Of course, the truncation of these $\delta_c$'s on $[0,\infty)$ would be more competitive with $\delta_1$ for $\theta \geq 0$, as well as with the maximum likelihood estimator $\max\{0,X\}$.
\end{example}

\begin{figure}[!ht]
  \centering
    \includegraphics[width=0.79\textwidth]{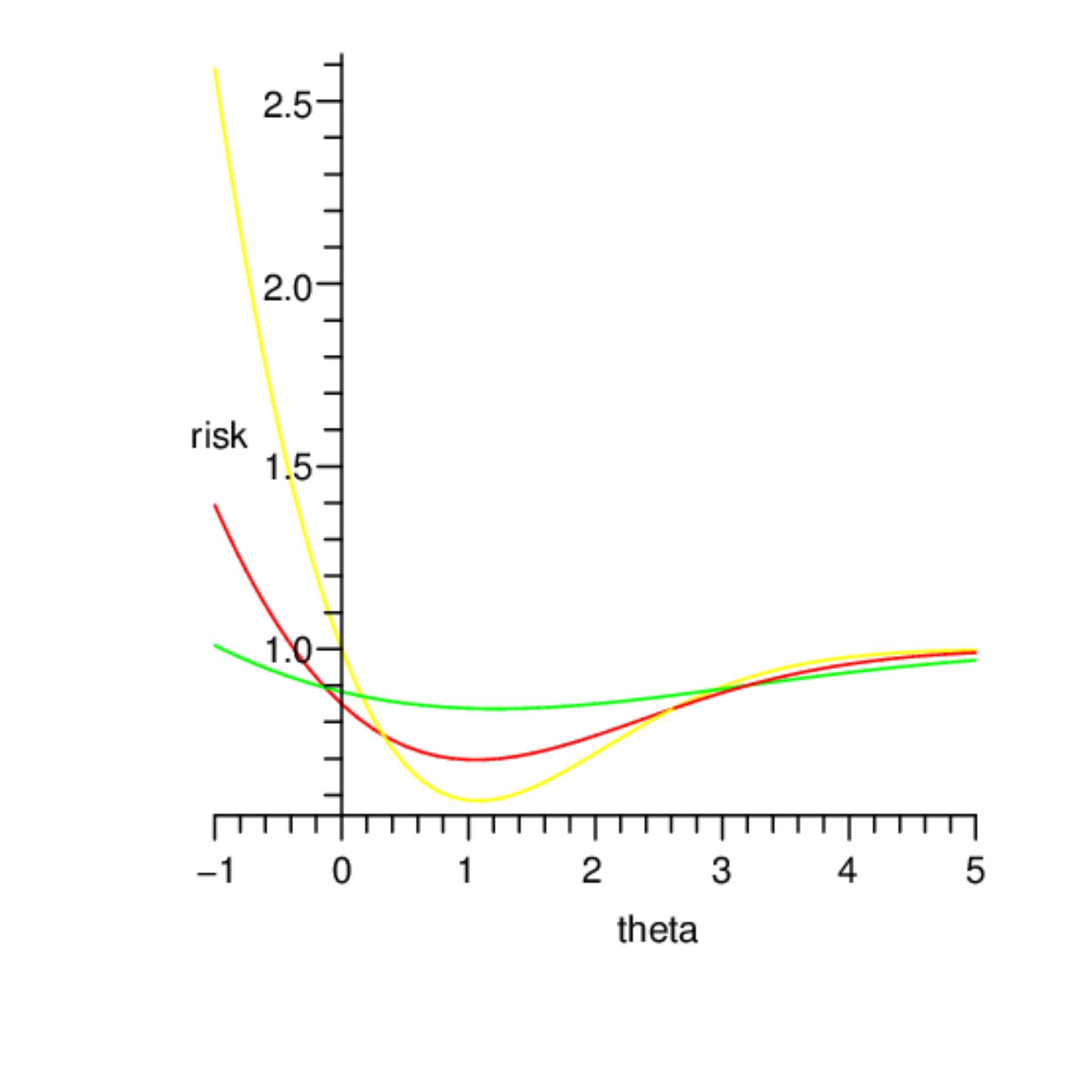}
     \label{ab2}
  \caption{Squared error loss as functions of $\theta$ for $\delta_c$, $c=1/2,3/4,1$}
\end{figure}

\section{Concluding Remarks}

We have presented a simple hierarchical model for Bayesian inference about a lower-bounded parameter $\theta$ with prior uncertainty on the lower bound $\alpha$, and with inference on $\alpha$ itself also considered.   For a normal model $X \sim N(\theta, \sigma^2)$ and the hierarchical prior $\pi(\theta|\alpha) = \mathbb{I}_{[\alpha,\infty)}(\theta)\,,\, \alpha \sim N(0,\sigma^2_{\alpha})$, we provided posterior analysis involving skewed normal distributions and showed that Bayes point estimators are minimax under squared error loss for the restricted parameter space $\theta \geq 0$.   Many of these features are appealing and extensions to doubly-bounded cases $\theta \in [m_1,m_2]$ with uncertain $m_1, m_2$, as well as multivariate extensions, for the normal and other models, merit further investigation.   It would be equally of interest to investigate in such settings, as we have done here focussing on minimaxity, the performance of some of the Bayes procedures from a frequentist perspective.

\section*{Acknowledgements}

Author Marchand gratefully acknowledges research support of the Natural Sciences and Engineering Research Council of Canada.   We are grateful to Latifa Ben Hadj Slimene for useful discussions and related numerical evaluations.

\end{document}